\documentclass{amsart2000}

\usepackage{xy,amssymb,
diagrams,times}
\usepackage{amsmath2000}
\xyoption{poly}
\xyoption{arc}
\xyoption{knot}
\xyoption{curve}
\xyoption{arrow}
\xyoption{all}

\newcommand{\eps}{\epsilon}

\newcommand{\C}{\mathbb{C}}

\newcommand{\Z}{\mathbb{Z}}

\newcommand{\Oscr}{{\mathcal O}}

\newcommand{\vtx}[1]{*+[o][F-]{\scriptscriptstyle #1}}
\newcommand{\cvtx}[1]{*+[F]{\scriptscriptstyle #1}}

\newcommand{\wis}[1]{{\text{\em \usefont{OT1}{cmtt}{m}{n} #1}}}

\newtheorem{definition}{Definition}
\newtheorem{theorem}{Theorem}
\newtheorem{proposition}{Proposition}
\newtheorem{lemma}{Lemma}

\title{Smooth Order Singularities}

\author{Raf Bocklandt, Lieven Le Bruyn and Geert Van de Weyer}
\address{Universiteit Antwerpen,
Antwerp (Belgium)}

\begin{document}

\maketitle

\begin{abstract}
In \cite{LBsmooth} it was shown that the center of Cayley-Hamilton smooth orders is smooth whenever the central dimension is at most two and that there may be singularities in higher dimensions. In this paper, we give methods to classify  central singularities of smooth orders up to smooth equivalence in arbitrary dimension and show that these methods are strong enough to complete the classification in dimension $\leq 6$. In particular we show that there is exactly one possible singularity type in dimension three : the conifold singularity. In dimensions $4$ (resp. $5$,$6$) there are precisely
$3$ (resp. $10$,$53$) types of singularities. This version of the paper contains the general techniques and the classification in dimension $\leq 4$. The full version (28 pages) containing the classifications in dimensions 5 and 6 is available from {\tt ftp://wins.uia.ac.be/pub/preprints/02/SOSfull.pdf}
\end{abstract}

\section{Introduction}  \label{intro}

One can define smoothness for a noncommutative algebra either by extending the homological
(Serre) or the categorical (Grothendieck) characterization of commutative regular algebras to the noncommutative world. In this paper we follow the second approach, started off by W. Schelter
\cite{Schelter} and C. Procesi \cite{ProcesiCH}, as we have an \'etale local description of these
{\em Cayley-Hamilton smooth orders} by the results of \cite{LBsmooth}. This local structure then
gives restrictions on the central simple algebras possessing a noncommutative smooth model.

An {\em algebra with trace map} $(A,tr)$ is an associative $\C$-algebra having a linear trace map
$tr~:~A \rTo A$ satisfying $tr(ab)=tr(ba)$, $tr(a)b=btr(a)$ and $tr(tr(a)b)=tr(a)tr(b)$. Morphisms in the category of algebras with trace are trace preserving $\C$-algebra morphisms. One has the identity
\[
\prod_{i=1}^n (t-x_i) = \sum_{i=0}^n (-1)^i \sigma_i t^{n-i} \]
where the $\sigma_i$ are the elementary symmetric polynomials in the $x_i$. There is another
generating set of the symmetric polynomials given by the power sums $\tau_k = \sum_i x_i^k$, so there are polynomials with rational coefficients
$
\sigma_k = p_k(\tau_1,\tau_2,\hdots,\tau_n)$ and we define the function $\sigma_k$ formally on any algebra with trace $(A,tr)$ to be
\[
\sigma_k(a) = p_k(tr(a),tr(a^2),\hdots,tr(a^n)) \]
This allows us to define a {\em formal $n$-th Cayley-Hamilton polynomial} for $(A,tr)$ by
\[
\chi_{n,a}(t) = \sum_{i=0}^n (-1)^i \sigma_i(a)t^{n-i} \]
and we say that $(A,tr)$ is an $n$-th Cayley-Hamilton algebra (or that $A \in \wis{alg@n}$) if
\[
tr(1) = n \qquad \text{and} \qquad \chi_{n,a}(a) = 0~\text{in $A$ for all $a \in A$} \]
The archetypical example of an $n$-th Cayley-Hamilton algebra is an order over a normal domain in a central simple algebra of degree $n$.

A {\em Cayley-Hamilton smooth algebra} is an affine $\C$-algebra in $\wis{alg@n}$ satisfying
Grothendieck's lifting characterization with respect to test-objects $(B,I)$ in $\wis{alg@n}$, that is,
any trace preserving algebra map $\phi$
\[
\begin{diagram}
A & \rDotsto^{\tilde{\phi}} & B \\
& \rdTo_{\phi} & \dOnto \\
& & \tfrac{B}{I}
\end{diagram}
\]
can be lifted to a trace preserving algebra map $\tilde{\phi}$ completing the diagram. C. Procesi
proved in \cite{ProcesiCH} that this categorical condition is equivalent to the geometric statement that the scheme $\wis{trep}_n~A$ of trace preserving $n$-dimensional representations of $A$ is a
smooth affine variety (though it may have several connected components). Moreover, the algebraic
quotient variety
\[
\wis{tiss}_n~A = \wis{trep}_n~A // GL_n \]
with respect to the natural base-change action has as its coordinate ring the central subalgebra
$tr(A)$ and its geometric points parametrize the trace preserving semi-simple $n$-dimensional representations of $A$. Of particular interest to us is the case of Cayley-Hamilton smooth {\em orders}, that is, when there is a Zariski open subset of $\wis{tiss}_n~A$ corresponding to simple
$n$-dimensional representations and (consequently) that $tr(A) = Z(A)$ the center of $A$.

If $A$ is a Cayley-Hamilton smooth order and $m$ is a maximal central ideal, then one can use the Luna slice theorem to determine the algebra structures of the $m$-adic completions (the \'etale local structure)
\[
\hat{A}_m \qquad \text{and} \qquad \widehat{Z(A)}_m \]
in terms of a {\em marked} quiver setting $(Q^{\bullet},\alpha)$, see \cite{LBsmooth}. To be precise, let $m$ be the point of $\wis{tiss}_n~A$ corresponding to the trace preserving semi-simple $n$-dimensional representation
\[
M = S_1^{\oplus e_1} \oplus \hdots \oplus S_k^{\oplus e_k} \]
where the $S_i$ are simple $d_i$-dimensional representations of $A$ occurring with multiplicity $e_i$ whence $n = \sum e_id_i$. Consider the quiver $Q$ on $k$ vertices $\{ v_1,\hdots, v_k \}$
(corresponding to the distinct simple components) such that the number of directed arrows from
$v_i$ to $v_j$ is given by the dimension of the subspace of the extension space $Ext^1_A(S_i,S_j)$ consisting of trace preserving algebra maps, see \cite{LBsmooth}. Let $\alpha$
be the dimension vector $(e_1,\hdots,e_k)$ (given by the multiplicities), then the $GL_n$-\'etale structure of $\wis{trep}_n~A$ in a neighborhood of the orbit $\Oscr(M)$ is isomorphic to the associated fiber bundle
\[
GL_n \times^{GL(\alpha)} \wis{rep}_{\alpha}~Q^{\bullet} \]
where $GL(\alpha) \rInto GL_n$ is determined by the dimensions $d_i$ and where $\wis{rep}_{\alpha}~Q^{\bullet}$ is the vectorspace of all $\alpha$-dimensional representations of the marked quiver $Q^{\bullet}$ (this means that some of the loops in $Q$ acquire a marking imposed by the trace preserving linear conditions, a representation of $Q^{\bullet}$ is a representation of $Q$ such that the matrix corresponding to a marked loop has trace zero). In particular, this implies that $\hat{A}_m$ is Morita equivalent to the completion of the algebra of
$GL(\alpha)$-equivariant maps $\wis{rep}_{\alpha}~Q^{\bullet} \rTo M_n(\C)$ at the maximal ideal corresponding to the zero representation and that $\hat{Z(A)}_m$ is isomorphic to the completion
\[
\C[[\wis{rep}_{\alpha}~Q^{\bullet}]]^{GL(\alpha)} \]
of the ring of polynomial quiver invariants at the maximal graded ideal. This fact allows us to study the central singularities of Cayley-Hamilton smooth orders. In \cite{LBsmooth} it was shown that the
center is smooth whenever the Krull dimension of the smooth order is $\leq 2$ and that there are central singularities possible in dimensions $\geq 3$.

Recall that two commutative local rings $C_m$ and $D_n$ are said to be {\em smooth equivalent}
if there are numbers $k$ and $l$ such that
\[
\hat{C}_m [[x_1,\hdots,x_k]] \simeq \hat{D}_n [[y_1,\hdots,y_l]] \]
A classification of all commutative singularities up to smooth equivalence is a hopeless task. Still,
because central singularities of Cayley-Hamilton smooth orders are determined by quiver invariants we will prove methods to attack this classification problem in principle and illustrate the
methods by giving a full classification in dimensions $\leq 6$. The main result of this paper is

\begin{theorem} Let $d$ be the dimension of the central variety $\wis{tiss}_n~A$ of a Cayley-Hamilton smooth order $A$. Then, if $d \leq 2$, $\wis{tiss}_n~A$ is smooth. If $d=3$ (resp. $4,5,6$) there are
exactly one (resp. three, ten and fifty three) types of central singularities possible.
\end{theorem}

In dimension three, the only possible central singularity is the so called {\em conifold singularity}
\[
\C[[u,v,x,y]]/(uv-xy) \]
In section two we give a general strategy to classify smooth equivalence classes of central singularities in any dimension, based on the reduction steps of \cite{Bocklandt} in the classification of the smooth quiver settings. In section three we give the proofs of the claims made and in the final two sections we give the details of the remaining classification result in dimensions $5$ and $6$.

\section{The strategy}

By the \'etale classification it suffices to classify marked quiver settings up to smooth equivalence,
that is, we want to determine when
\[
\C[\wis{rep}_{\alpha_1}~Q_1^{\bullet}]^{GL(\alpha_1)}[x_1,\hdots,x_k] \simeq
\C[\wis{rep}_{\alpha_2}~Q_2^{\bullet}]^{GL(\alpha_2)}[y_1,\hdots,y_l]
\]
In the case of quivers, a full classification of all the quiver settings $(Q,\alpha)$ such that the ring of invariants is a polynomial ring was given in \cite{Bocklandt}. The proof relies on a number of reduction steps which modify the ring of invariants only up to polynomial extensions. We will recall these reduction steps as well as their obvious extensions to marked quivers. In the quiver diagrams below, the vertex-dimension component is depicted in the vertex and the number of multiple arrows between two vertices is given by a superscript, unless this number is $\leq 3$ in which case the number of arrows is drawn. In the diagrams below we only depict the quiver-neighborhood of the vertex where a change is made, the remaining part of the quiver setting is left unchanged.

Recall that the {\em Euler form} $\chi_Q$ of a quiver $Q$ is the bilinear form on $\Z^k$ (if $Q$ has $k$ vertices) determined by the integral $k \times k$ matrix having as its $(i,j)$-entry
\[
\chi_{Q,ij} = \delta_{ij} - \# \{ \text{arrows from $v_i$ to $v_j$} \} \]
With $\epsilon_v$ we denote the basevector concentrated in vertex $v$ and $\alpha_v$ will denote the vertex dimension component of $\alpha$ in vertex $v$. There are three types of reduction moves, each with their own condition and effect on the ring of invariants.

{\bf Vertex removal : } Let $(Q^{\bullet},\alpha)$ be a marked quiver setting and $v$ a vertex  satisfying the condition $C^v_V$, that is, $v$ is without (marked) loops and satisfies
\[
\chi_Q(\alpha,\epsilon_v) \geq 0 \quad \text{or} \quad \chi_Q(\epsilon_v,\alpha) \geq 0 \]
Define the new quiver setting $(Q^{\bullet'},\alpha')$ obtained by the operation $R^v_V$ which removes the vertex $v$ and composes all arrows through $v$, the dimensions of the other vertices are unchanged :
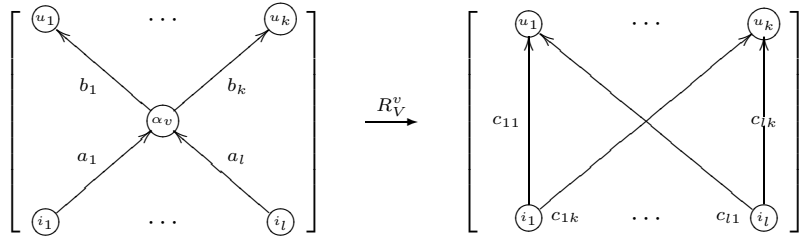
\[
\left[ ~\vcenter{
\xymatrix@=1cm{
\vtx{u_1}&\cdots &\vtx{u_k}\\
&\vtx{\alpha_v}\ar[ul]^{b_1}\ar[ur]_{b_k}&\\
\vtx{i_1}\ar[ur]^{a_1}&\cdots &\vtx{i_l}\ar[ul]_{a_l}}}
~\right] \quad
\rTo^{R^v_V} \quad
\left[~\vcenter{
\xymatrix@=1cm{
\vtx{u_1}&\cdots &\vtx{u_k}\\
&&\\
\vtx{i_1}\ar[uu]^{c_{11}}\ar[uurr]_<<{c_{1k}}&\cdots &\vtx{i_l}\ar[uu]|{c_{lk}}\ar[uull]^<<{c_{l1}}}}
~\right].
\]
where $c_{ij} = a_ib_j$ (observe that some of the incoming and outgoing vertices may be the
same so that one obtains loops in the corresponding vertex). In this case we have
\[
\C[\wis{rep}_{\alpha}~Q^{\bullet}]^{GL(\alpha)} \simeq \C[\wis{rep}_{\alpha'}~Q^{\bullet'}]^{GL(\alpha')} \]

{\bf loop removal :} Let $(Q^{\bullet},\alpha)$ be a marked quiver setting and $v$ a vertex satisfying the condition $C^v_l$  that the vertex-dimension $\alpha_v = 1$ and there are $k \geq 1$ loops in
$v$. Let $(Q^{\bullet'},\alpha)$ be the quiver setting obtained by the loop removal operation $R^v_l$
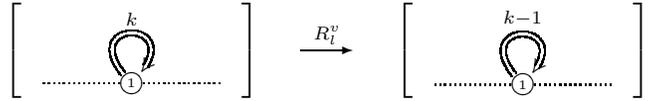
\[
\left[~\vcenter{
\xymatrix@=1cm{
&\vtx{1}\ar@{..}[r]\ar@{..}[l]\ar@(lu,ru)@{=>}^k&}}
~\right]\quad \rTo^{R^v_l} \quad
\left[~\vcenter{
\xymatrix@=1cm{
&\vtx{1}\ar@{..}[r]\ar@{..}[l]\ar@(lu,ru)@{=>}^{k-1}&}}
~\right].\]
removing one loop in $v$ and keeping the dimension vector the same, then
\[
\C[\wis{rep}_{\alpha}~Q^{\bullet}]^{GL(\alpha)} \simeq \C[\wis{rep}_{\alpha}~Q^{\bullet'}]^{GL(\alpha)}[x] \]

{\bf Loop removal : } Let $(Q^{\bullet},\alpha)$ be a marked quiver setting and $v$ a vertex satisfying condition $C^v_L$, that is, the vertex dimension $\alpha_v \geq 2$, $v$ has precisely one (marked) loop in $v$ and
\[
\chi_Q(\epsilon_v,\alpha) = -1 \quad \text{or} \quad \chi_Q(\alpha,\epsilon_v) = -1 \]
(that is, there is exactly one other incoming or outgoing arrow from/to  a vertex with dimension $1$). Let $(Q^{\bullet'},\alpha)$ be the marked quiver setting obtained by changing the quiver as
indicated below (depending on whether the incoming or outgoing condition is satisfied and whether there is a loop or a marked loop in $v$)
\vspace{.5cm}
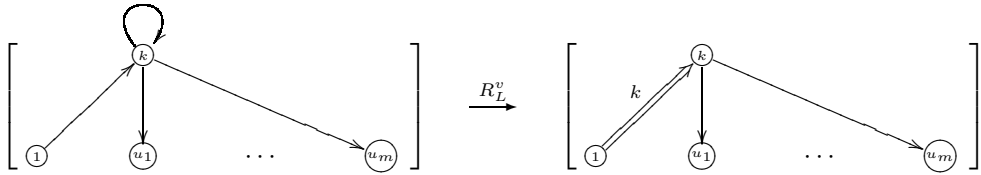
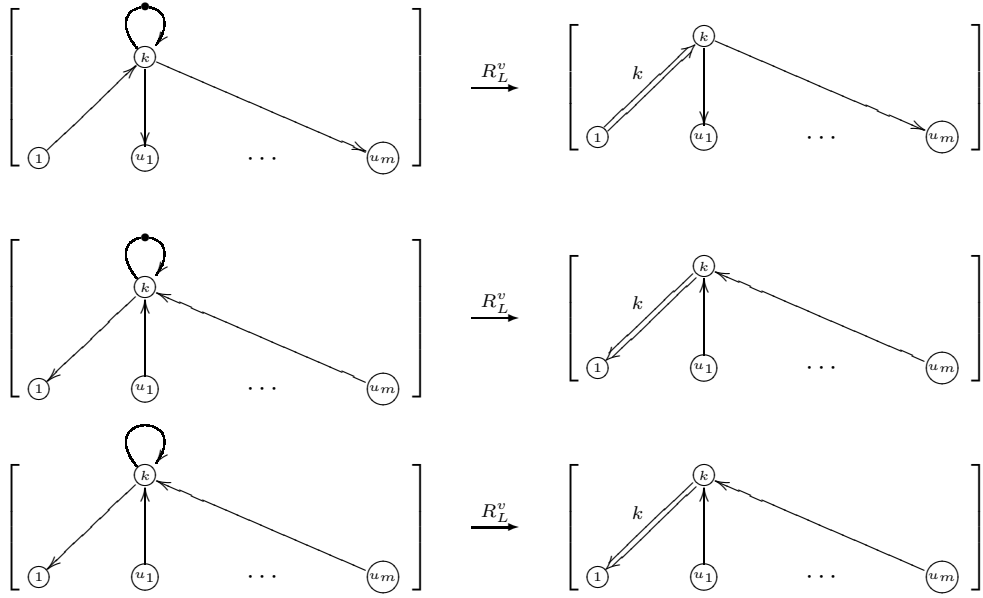
\[
\left[~\vcenter{
\xymatrix@=1cm{
&\vtx{k}\ar[d]\ar[drr]\ar@(lu,ru)&&\\
\vtx{1}\ar[ur]&\vtx{u_1}&\cdots &\vtx{u_m}}}
~\right]\quad \rTo^{R^v_L} \quad
\left[~\vcenter{
\xymatrix@=1cm{
&\vtx{k}\ar[d]\ar[drr]&&\\
\vtx{1}\ar@2[ur]^k&\vtx{u_1}&\cdots &\vtx{u_m}}}
~\right]
\]
\vspace{.5cm}
\[
\left[~\vcenter{
\xymatrix@=1cm{
&\vtx{k}\ar[d]\ar[drr]\ar@(lu,ru)|{\bullet}&&\\
\vtx{1}\ar[ur]&\vtx{u_1}&\cdots &\vtx{u_m}}}
~\right]\quad \rTo^{R^v_L} \quad
\left[~\vcenter{
\xymatrix@=1cm{
&\vtx{k}\ar[d]\ar[drr]&&\\
\vtx{1}\ar@2[ur]^k&\vtx{u_1}&\cdots &\vtx{u_m}}}
~\right] \]
\vspace{.5cm}
\[
\left[~\vcenter{
\xymatrix@=1cm{
&\vtx{k}\ar[dl]\ar@(lu,ru)|{\bullet}&&\\
\vtx{1}&\vtx{u_1}\ar[u]&\cdots &\vtx{u_m}\ar[ull]}}
~\right]\quad \rTo^{R^v_L} \quad
\left[~\vcenter{
\xymatrix@=1cm{
&\vtx{k}\ar@2[dl]_{k}&&\\
\vtx{1}&\vtx{u_1}\ar[u]&\cdots &\vtx{u_m}\ar[ull]}}
~\right].
\]
\vspace{.5cm}
\[
\left[~\vcenter{
\xymatrix@=1cm{
&\vtx{k}\ar[dl]\ar@(lu,ru)&&\\
\vtx{1}&\vtx{u_1}\ar[u]&\cdots &\vtx{u_m}\ar[ull]}}
~\right]\quad \rTo^{R^v_L} \quad
\left[~\vcenter{
\xymatrix@=1cm{
&\vtx{k}\ar@2[dl]_k&&\\
\vtx{1}&\vtx{u_1}\ar[u]&\cdots &\vtx{u_m}\ar[ull]}}
~\right].
\]

and the dimension vector is left unchanged, then we have
\[
\C[\wis{rep}_{\alpha}~Q^{\bullet}]^{GL(\alpha)} = \begin{cases}
\C[\wis{rep}_{\alpha}~Q^{\bullet'}]^{GL(\alpha)}[x_1,\hdots,x_k]~\quad \text{(loop)} \\
\C[\wis{rep}_{\alpha}~Q^{\bullet'}]^{GL(\alpha)}[x_1,\hdots,x_{k-1}]~\quad \text{(marked loop)} 
\end{cases}
\]

\begin{definition} A marked quiver $Q^{\bullet}$ is said to be strongly connected if for every pair of vertices $\{ v,w \}$ there is an oriented path from $v$ to $w$ and an oriented path from $w$ to $v$.

A marked quiver setting $(Q^{\bullet},\alpha)$ is said to be {\em reduced} if and only if there is {\em no} vertex $v$ such that one of the conditions $C^v_V$, $C^v_l$ or $C^v_L$ is satisfied.
\end{definition}

\begin{lemma} Every marked quiver setting $(Q^{\bullet}_1,\alpha_1)$ can be reduced by a sequence of operations $R^v_V,R^v_l$ and $R^v_L$ to a {\em reduced} quiver setting $(Q^{\bullet}_2,\alpha_2)$ such that 
\[
\C[\wis{rep}_{\alpha_1}~Q^{\bullet}_1]^{GL(\alpha_1)} \simeq
\C[\wis{rep}_{\alpha_2}~Q^{\bullet}_2]^{GL(\alpha_2)}[x_1,\hdots,x_z] \]
Moreover, the number $z$ of extra indeterminates is determined by the reduction sequence
\[
(Q^{\bullet}_2,\alpha_2) = R^{v_{i_u}}_{X_u} \circ \hdots \circ R^{v_{i_1}}_{X_1} (Q^{\bullet}_1,\alpha_1) \]
where for every $1 \leq j \leq u$, $X_j \in \{ V,l,L \}$. More precisely,
\[
z = \sum_{X_j = l} 1 + \sum_{X_j = L}^{(unmarked)} \alpha_{v_{i_j}} + \sum_{X_j = L}^{(marked)} (\alpha_{v_{i_j}} - 1) \]
\end{lemma}

\begin{proof}
As any reduction step removes a (marked) loop or a vertex, any sequence of reduction steps starting with $(Q^{\bullet}_1,\alpha_1)$ must eventually end in a reduced marked quiver setting. The
statement then follows from the discussion above.
\end{proof}

\par \vskip 3mm
As the reduction steps have no uniquely determined inverse, there is no a priori reason why the
reduced quiver setting of the previous lemma should be unique. Nevertheless this is true as we will prove in section~\ref{unique} :

\begin{theorem} \label{thmunique} Every marked quiver setting $(Q^{\bullet}_1,\alpha_1)$ can be transformed by a sequence of reduction steps $R^v_V,R^v_l$ or $R^v_L$ to a {\em uniquely determined} reduced
marked quiver setting $(Q^{\bullet}_2,\alpha_2)$.
\end{theorem}

This result shows that it is enough to classify {\em reduced} marked quiver settings up to smooth equivalence. We can always assume that the quiver $Q$ is strongly connected (if not, the ring of invariants is the tensor product of the rings of invariants of the maximal strongly connected subquivers). Our aim is to classify the reduced quiver singularities up to equivalence, so we need to determine the Krull dimension of the rings of invariants.

\begin{lemma} Let $(Q^{\bullet},\alpha)$ be a reduced marked quiver setting and $Q$ strongly connected. Then,
\[
dim~\wis{iss}_{\alpha}~Q^{\bullet} = 1 - \chi_Q(\alpha,\alpha) - m \]
where $m$ is the total number of marked loops in $Q^{\bullet}$.
\end{lemma}

\begin{proof}
Because $(Q^{\bullet},\alpha)$ is reduced, none of the vertices satisfies condition $C^v_V$, whence
\[
\chi_Q(\epsilon_v,\alpha) \leq -1 \quad \text{and} \quad \chi_Q(\alpha,\epsilon_v) \leq -1 \]
for all vertices $v$. In particular it follows (because $Q$ is strongly connected) from \cite{LBProcesi}
that $\alpha$ is the dimension vector of a simple representation of $Q$ and that the dimension of the quotient variety
\[
dim~\wis{iss}_{\alpha}~Q = 1 - \chi_Q(\alpha,\alpha) \]
Finally, separating traces of the loops to be marked gives the required formula.
\end{proof}

\par \vskip 3mm
Applying the main result of \cite{Bocklandt} we have all marked quiver settings having a regular ring of invariants. This result also describes the smooth locus of the central variety of a Cayley-Hamilton smooth order using the \'etale local description of section~\ref{intro}.

\begin{theorem} Let $(Q^{\bullet},\alpha)$ be a marked quiver setting such that $Q$ is strongly connected. Then $\wis{iss}_{\alpha}~Q^{\bullet}$ is smooth if and only if the unique reduced marked quiver setting to which $(Q^{\bullet},\alpha)$ can be reduced is one of the following five types
\[
\xymatrix@=1cm{
\vtx{k} & & \vtx{k}\ar@(ul,dl) & & \vtx{2}\ar@(ul,dl)\ar@(ur,dr) &&  \vtx{2}\ar@(ul,dl)\ar@(ur,dr)|{\bullet} &&
\vtx{2}\ar@(ul,dl)|{\bullet}\ar@(ur,dr)|{\bullet}}
\]
\end{theorem}

\begin{proof}
Because the ring of invariants is graded it suffices to prove smoothness in the origin. Consider the underlying quiver $Q$, apply the main result of \cite{Bocklandt} and separate traces of the marked loops.
\end{proof}

\par \vskip 3mm
The next step is to classify for a given dimension $d$ all reduced marked quiver settings
$(Q^{\bullet},\alpha)$ such that $dim~\wis{iss}_{\alpha}~Q^{\bullet} = d$. The following result
limits the possible cases drastically in low dimensions.

\begin{lemma} \label{counting} Let $(Q^{\bullet},\alpha)$ be a reduced marked quiver setting on
$k \geq 2$ vertices. Then,
\[
dim~\wis{iss}_{\alpha}~Q^{\bullet} \geq 1 + \sum_{\xymatrix@=1cm{ \vtx{a} }}^{a \geq 1} a + 
\sum_{\xymatrix@=1cm{ \vtx{a}\ar@(ul,dl)|{\bullet} }}^{a > 1}(2a-1) +
 \sum_{\xymatrix@=1cm{ \vtx{a}\ar@(ul,dl)}}^{a > 1}(2a) + \sum_{\xymatrix@=1cm{ \vtx{a}\ar@(ul,dl)|{\bullet}\ar@(ur,dr)|{\bullet}}}^{a > 1} (a^2+a-2) + \]
\[
\sum_{\xymatrix@=1cm{ \vtx{a}\ar@(ul,dl)|{\bullet}\ar@(ur,dr)}}^{a > 1} (a^2+a-1) +
\sum_{\xymatrix@=1cm{ \vtx{a}\ar@(ul,dl)\ar@(ur,dr)}}^{a > 1} (a^2+a) + \hdots +
\sum_{\xymatrix@=1cm{ \vtx{a}\ar@(ul,dl)|{\bullet}_{k}\ar@(ur,dr)^{l}}}^{a > 1} ((k+l-1)a^2+a-k) + \hdots
\]
In this sum the contribution of a vertex $v$ with $\alpha_v = a$ is determined by the number of
(marked) loops in $v$. By the reduction steps (marked) loops only occur at vertices where
$\alpha_v > 1$.
\end{lemma}

\begin{proof}
We know that the dimension of $\wis{iss}_{\alpha}~Q^{\bullet}$ is equal to
\[
1 - \chi_Q(\alpha,\alpha) - m = 1 - \sum_v \chi_Q(\epsilon_v,\alpha) \alpha_v - m \]
If there are no (marked) loops at $v$, then $\chi_Q(\epsilon_v,\alpha) \leq -1$ (if not we would reduce further)
which explains the first sum. If there is exactly one (marked) loop at $v$ then $\chi_Q(\epsilon_v,\alpha) \leq -2$ for if $\chi_Q(\epsilon_v,\alpha) = -1$ then there is just one
outgoing arrow to a vertex $w$ with $\alpha_w = 1$ but then we can reduce the quiver setting
further. This explains the second and third sums. If there are $k$ marked loops and $l$ ordinary
loops in $v$ (and $Q$ has at least two vertices) , then
\[
-\chi_Q(\epsilon_v,\alpha)\alpha_v - k \geq ((k+l)\alpha_v -\alpha_v + 1)\alpha_v - k \]
which explains all other sums.
\end{proof}

\par \vskip 3mm
Observe that  the dimension of the quotient variety of the one vertex marked quivers
\[
\xymatrix@=1cm{ \vtx{a}\ar@(ul,dl)|{\bullet}_{k}\ar@(ur,dr)^{l}} \]
is equal to $(k+l-1)a^2+1-k$ and is singular (for $a \geq 2$) unless $k+l = 2$. We will now classify
the reduced singular settings when there are at least two vertices in low dimensions. By the previous lemma it follows immediately that
\begin{enumerate}
\item{the maximal number of vertices in a reduced marked quiver setting $(Q^{\bullet},\alpha)$ of
dimension $d$ is $d-1$ (in which case all vertex dimensions must be equal to one)}
\item{if a vertex dimension in a reduced marked quiver setting is $a \geq 2$, then the dimension
$d \geq 2a$.}
\end{enumerate}

\begin{lemma} \label{prop1} Let $(Q^{\bullet},\alpha)$ be a reduced marked quiver setting such that
$\wis{iss}_{\alpha}~Q^{\bullet}$ is singular of dimension $d \leq 5$, then
$\alpha = (1,\hdots,1)$. Moreover, each vertex must have at least two incoming and two outgoing
arrows and no loops.
\end{lemma}

\begin{proof}
From the lower bound of the sum formula it follows that if some $\alpha_v > 1$ it must be
equal to $2$ and must have a unique marked loop and there can only be one other vertex $w$
with $\alpha_w = 1$. If there are $x$ arrows from $w$ to $v$ and $y$ arrows from $v$ to $w$,
then
\[
dim~\wis{iss}_{\alpha}~Q^{\bullet} = 2(x+y)-1 \]
whence $x$ or $y$ must be equal to $1$ contradicting reducedness. The second statement follows
as otherwise we could perform extra reductions.
\end{proof}

\begin{proposition}
The only reduced marked quiver singularity in dimension 3 is
\[
3_{con}~:~\vcenter{ \xymatrix{\vtx{1}\ar@2@/^/[rr]&&\vtx{1}\ar@2@/^/[ll]} }
\]
The reduced marked quiver singularities in dimension 4 are
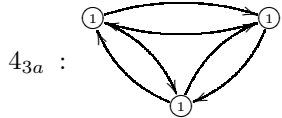
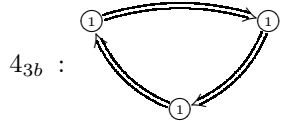
\[
4_{3a}~:~\vcenter{
\xymatrix{\vtx{1}\ar@/^/[rr]\ar@/^/[rd]&&\vtx{1}\ar@/^/[ll]\ar@/^/[ld]\\
&\vtx{1}\ar@/^/[ru]\ar@/^/[lu]&} } \qquad 4_{3b}~:~ \vcenter{
\xymatrix{\vtx{1}\ar@2@/^/[rr]&&\vtx{1}\ar@2@/^/[ld]\\
&\vtx{1}\ar@2@/^/[lu]&} }  \qquad 4_{2}~:~ \vcenter{
\xymatrix{\vtx{1}\ar@2@/^/[rr]&&\vtx{1}\ar@3@/^/[ll]} }
\]
\end{proposition}

\begin{proof}
All one vertex marked quiver settings with quotient dimension 
$\leq 5$ are smooth, so we are in the situation of  lemma~\ref{prop1}.
If the dimension is $3$ there must be two vertices each having exactly two incoming and
two outgoing arrows, whence the indicated type is the only one. The resulting singularity
is the {\em conifold singularity}
\[
\frac{\C[[x,y,u,v]]}{(xy-uv)} \]
In dimension $4$ we can have three or two vertices. In the first case, each vertex must have exactly
two incoming and two outgoing arrows whence the first two cases. If there are two vertices, then
just one of them has three incoming arrows and one has three outgoing arrows.
\end{proof}

\par \vskip 3mm
In dimensions $5$ and $6$ one can give a classification of all reduced singularities by hand (see
the full version of this paper). This concludes the first step in our strategy, the next will be to
distinguish reduced singularities of the same dimension up to (\'etale) isomorphism.

\section{Fingerprinting singularities}  \label{finger}

In this section we will outline methods to distinguish two reduced marked quiver settings
$(Q^{\bullet}_1,\alpha_1)$ and $(Q_2^{\bullet},\alpha_2)$ having the same quotient dimension $d$. Recall from \cite{LBProcesi} that the rings of quiver invariants are generated by taking traces along oriented cycles in the quiver (again separating traces gives the same result for marked quivers). Assume that all vertex dimensions are equal to one, then one can write any (trace of an)
oriented cycle as a product of (traces of) {\em primitive} oriented cycles (that is, those that cannot
be decomposed further). From this one deduces immediately :

\begin{lemma} Let $(Q^{\bullet},\alpha)$ be a reduced marked quiver setting such that all $\alpha_v = 1$.
Let $m$ be the maximal graded ideal of $\C[\wis{rep}_{\alpha}~Q^{\bullet}]^{GL(\alpha)}$, then a
vectorspace basis of 
\[
\tfrac{m^i}{m^{i+1}} \]
is given by the oriented cycles in $Q$ which can be written as a product of $i$ primitive cycles but
not as a product of $i+1$ such cycles.
\end{lemma}

Clearly, the dimensions of the quotients $m^i/m^{i+1}$ are (\'etale) isomorphism invariants. Hence, for
$d \leq 5$ this simple minded counting method can be used to separate quiver singularities.

\begin{theorem} There are precisely three reduced quiver singularities in dimension $d=4$.
\end{theorem}

\begin{proof}
The number of primitive oriented cycles of the three types of reduced marked quiver settings in dimension four
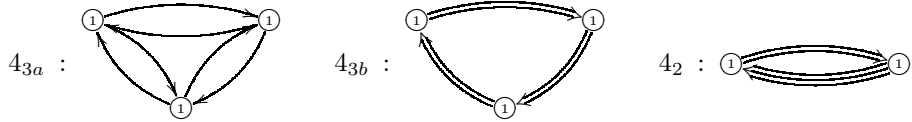
\[
4_{3a}~:~\vcenter{
\xymatrix{\vtx{1}\ar@/^/[rr]\ar@/^/[rd]&&\vtx{1}\ar@/^/[ll]\ar@/^/[ld]\\
&\vtx{1}\ar@/^/[ru]\ar@/^/[lu]&} } \qquad 4_{3b}~:~ \vcenter{
\xymatrix{\vtx{1}\ar@2@/^/[rr]&&\vtx{1}\ar@2@/^/[ld]\\
&\vtx{1}\ar@2@/^/[lu]&} }  \qquad 4_{2}~:~ \vcenter{
\xymatrix{\vtx{1}\ar@2@/^/[rr]&&\vtx{1}\ar@3@/^/[ll]} }
\]
is $5$, respectively $8$ and $6$. Hence, they give nonisomorphic rings of invariants.
\end{proof}

\par \vskip 3mm
In section~5 of the full version, we will classify the reduced quiver singularities for $d=5$. If some of the vertex dimensions are $\geq 2$ we have no easy description of the vectorspaces $m^i/m^{i+1}$
and we need a more refined argument. The idea is to answer the question "what other singularities can the reduced singularity see ?" by the theory of local quivers of \cite{LBProcesi}.

Let $Q$ be a quiver (we will indicate the necessary changes to be made for marked quivers below) and $\alpha$ a dimension vector. An $\alpha$-{\em representation type} is  a datum
\[
\tau = (e_1,\beta_1;\hdots;e_l,\beta_l) \]
where the $e_i$ are natural numbers $\geq 1$, the $\beta_i$ are dimension vectors of simple
representations of $Q$ (for which we have a precise description by \cite{LBProcesi}) such that
$\alpha = \sum_i e_i \beta_i$. Any neighborhood of the trivial representation contains semi-simple
representations of $Q$ of type $\tau$ for any $\alpha$-representation type. 

To determine the dimension of the corresponding strata and the nature of their singularities we construct a new quiver
$Q_{\tau}$, the {\em local quiver}, on $l$ vertices (the number of distinct simple components)
say $\{ w_1,\hdots,w_l \}$ such that the number of oriented arrows (or loops) from $w_i$ to $w_j$ is given by the number
\[
\delta_{ij} - \chi_Q(\beta_i,\beta_j) \]
There is an \'etale local isomorphism between a neighborhood of a semi-simple $\alpha$-dimensional representation of type $\tau$ and a neighborhood of the trivial representation of $\wis{iss}_{\alpha_{\tau}}~Q_{\tau}$ where
$\alpha_{\tau} = (e_1,\hdots,e_l)$ is the dimension vector determined by the multiplicities.

As a consequence we see that the dimension of the corresponding strata is equal to the number of
loops in $Q_{\tau}$. Now, assume that $\wis{iss}_{\alpha_{\tau}}~Q_{\tau}$ has a singularity, then the couple 
\[
(\text{dimension of strata, type of singularity})
\]
 is a characteristic feature of the singularity of $\wis{iss}_{\alpha}~Q$ and one can often distinguish types by these couples. In the case of a marked quiver one proceeds as before for the underlying quiver and in the final result compensates for the markings (that is, one marks as many loops in the local quiver in the vertices giving a non-zero contribution to the original marked vertex).

Recall from \cite{LBProcesi} that there is a partial ordering $\tau < \tau'$ on the $\alpha$-representation types induced by degeneration of representations. The {\em fingerprint}
of a reduced quiver singularity will be the Hasse diagram of those $\alpha$-representation types $\tau$ such that the local marked quiver setting $(Q^{\bullet}_{\tau},\alpha_{\tau})$ can be reduced
to a reduced quiver singularity (necessarily occurring in lower dimension and the difference between the two dimensions gives the dimension of the stratum).

Clearly, this method fails in case the marked quiver singularity is an {\em isolated singularity}.
Fortunately, we have a complete classification of such singularities by the work of \cite{BockSym}.

\begin{theorem}~\cite{BockSym}
The only reduced marked quiver settings $(Q^{\bullet},\alpha)$ such that the quotient variety is an isolated singularity are of the form
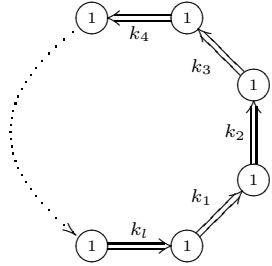
\[
\xy 0;/r.15pc/:
\POS (0,0) *+{\txt{\tiny $1$}}*\cir<6pt>{} ="a",
(20,0) *+{\txt{\tiny $1$}}*\cir<6pt>{} ="b",
(34,14) *+{\txt{\tiny $1$}}*\cir<6pt>{} ="c",
(34,34) *+{\txt{\tiny $1$}}*\cir<6pt>{} ="d",
(20,48) *+{\txt{\tiny $1$}}*\cir<6pt>{} ="e",
(0,48) *+{\txt{\tiny $1$}}*\cir<6pt>{} ="f"
\POS"a" \ar@{=>}^{k_l} "b"
\POS"b" \ar@{=>}^{k_1} "c"
\POS"c" \ar@{=>}^{k_2} "d"
\POS"d" \ar@{=>}^{k_3} "e"
\POS"e" \ar@{=>}^{k_4} "f"
\POS"f" \ar@/_7ex/@{.>} "a"
\endxy 
\]
where $Q$ has $l$ vertices and all $k_i \geq 2$. The dimension of the corresponding quotient
is
\[
d = \sum_i k_i + l - 1 \]
and the {\em unordered $l$-tuple} $\{ k_1,\hdots,k_l \}$ is an (\'etale) isomorphism invariant of the ring of invariants.
\end{theorem}

Not only does this result distinguish among isolated reduced quiver singularities, but it also shows that in all other marked quiver settings we will have additional families of singularities. We will illustrate the method in some detail to separate the reduced marked quiver settings in dimension $6$ having one vertex of dimension two.

\begin{proposition} \label{classif62} The reduced singularities of dimension $6$ such that $\alpha$ contains a component equal to $2$ are pairwise non-equivalent.
\end{proposition}

\begin{proof}
In section~6 of the full version, we will show that the relevant reduced marked quiver setting are the following
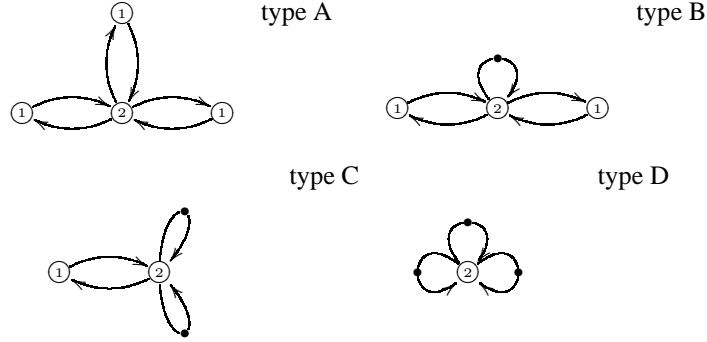
\[
\xymatrix@=1cm{&\vtx{1}\ar@/^/[d]& \\ \vtx{1}\ar@/^/[r]&\vtx{2}\ar@/^/[l]\ar@/^/[r]\ar@/^/[u]&\vtx{1}\ar@/^/[l]} \quad \text{type A} \qquad
\xymatrix@=1cm{& & \\ \vtx{1}\ar@/^/[r]&\vtx{2}\ar@/^/[l]\ar@/^/[r]\ar@(ul,ur)[]|{\bullet}&\vtx{1}\ar@/^/[l]}
\quad \text{type B} \]
\[
\xymatrix@=1cm{& \\ \vtx{1}\ar@/^/[r]&\vtx{2}\ar@/^/[l]\ar@(u,ur)[]|{\bullet}\ar@(d,dr)[]|{\bullet}&} \quad
\text{type C} \qquad
\xymatrix@=1cm{& \\ \vtx{2}\ar@(ur,dr)[]|{\bullet}\ar@(ul,dl)[]|{\bullet}\ar@(ul,ur)[]|{\bullet}&} \quad \text{type D}
\]
We will order the vertices such that $\alpha_1 = 2$.

\par \vskip 2mm
\noindent
{\bf type A : } There are three different 
representation types $\tau_1 = (1,(2;1,1,0);1,(0;0,0,1))$ (and permutations of the $1$-vertices).
The local quiver setting  has the form

\par \vskip 2mm
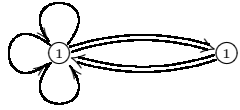
\[
\xymatrix{\vtx{1}\ar@2@/^/[rr]\ar@(ur,ul)\ar@(ul,dl)\ar@(dl,dr)&&\vtx{1}\ar@2@/^/[ll]}
\]

\par \vskip 4mm \noindent
because for $\beta_1=(2;1,1,0)$ and $\beta_2=(0;0,0,1)$ we have that $\chi_Q(\beta_1,\beta_1)=-2$, $\chi_Q(\beta_1,\beta_2)=-2$,
$\chi_Q(\beta_2,\beta_1)=-2$ and $\chi(\beta_2,\beta_2)=1$. These three representation types each give a three dimensional family of conifold (type $3_{con}$) singularities.

Further, there are three different representation types $
\tau_2 = (1,(1;1,1,0);1,(1;0,0,1))$
(and permutations) of which the 
local quiver setting is of the form

\par
\vskip 2mm
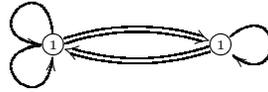
\[
\xymatrix{\vtx{1}\ar@2@/^/[rr]\ar@(u,l)\ar@(l,d)&&\vtx{1}\ar@2@/^/[ll]\ar@(ur,dr)}
\]

\par
\vskip 4mm \noindent
as with $\beta_1 = (1;1,1,0)$ and $\beta_2=(1;0,0,1)$ we have 
$\chi_Q(\beta_1,\beta_1) = -1$, $\chi_Q(\beta_1,\beta_2) = -2$, $\chi_Q(\beta_2,\beta_1) = -2$ and
$\chi_Q(\beta_2,\beta_2) = 0$. These three representation types each give a
three dimensional family of conifold singularities.

Finally, there are the three representation types 
\[
\tau_3 = (1,(1;1,0,0);1,(1;0,1,0);1,(0;0,0,1))
\]
(and permutations) with local quiver setting

\par \vskip 2mm
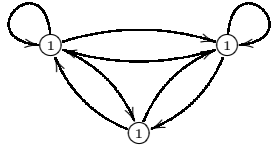
\[
\xymatrix{\vtx{1}\ar@/^/[rr]\ar@/^/[rd]\ar@(u,l)&&\vtx{1}\ar@/^/[ll]\ar@/^/[ld]\ar@(u,r)\\
&\vtx{1}\ar@/^/[ru]\ar@/^/[lu]&}
\]
These three types each give a two dimensional family of reduced singularities of type $4_{3a}$.

The degeneration order on representation types gives $\tau_1 < \tau_3$ and $\tau_2 < \tau_3$
(but for different permutations) and the {\em fingerprint} of this reduced singularity can be depicted as
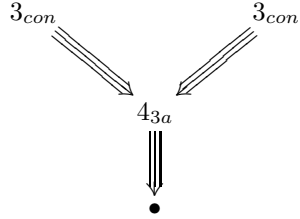
\[
\xymatrix{\txt{$3_{con}$}\ar@3[rd]& & \txt{$3_{con}$}\ar@3[ld] \\
& \txt{$4_{3a}$}\ar@3[d] & \\
& \bullet & } \]

\par \vskip 2mm
\noindent
{\bf type B : } There is one
representation type $\tau_1 = (1,(1;1,0);1,(1;0,1))$ giving as above a three dimensional family of conifold singularities, one representation type $\tau_2 = (1,(1;1,1);1,(1;0,0))$ giving a three dimensional family of conifolds and finally one representation type 
\[
\tau_3 = (1,(1;0,0);1,(1;0,0);1,(0;1,1);1,(0;0,1))
\]
of which the local quiver setting has the form
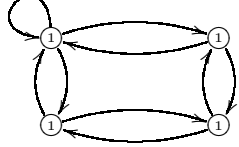
\[
\xymatrix{\vtx{1}\ar@/^/[rr]\ar@/^/[d]\ar@(u,l)&&\vtx{1}\ar@/^/[ll]\ar@/^/[d]\\
\vtx{1}\ar@/^/[rr]\ar@/^/[u]&&\vtx{1}\ar@/^/[ll]\ar@/^/[u]}
\]
(the loop in the downright corner is removed to compensate for the marking)
 giving rise to a one-dimensional family of five-dimensional singularities of type $5_{4a}$. This gives the fingerprint
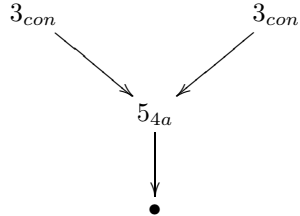
\[
\xymatrix{\txt{$3_{con}$}\ar[rd]& & \txt{$3_{con}$}\ar[ld] \\
& \txt{$5_{4a}$}\ar[d] & \\
& \bullet & } \]

\par \vskip 2mm
\noindent
{\bf type C : } We have a three dimensional family of conifold singularities coming from the representation type $(1,(1;1);1,(1;0))$ and a two-dimensional family of type $4_{3a}$ singularities corresponding to the representation type $(1,(1;0);1,(1,0);1,(0;1))$. Therefore, the
fingerprint is depicted as
\[
3_{con} \rTo 4_{3a} \rTo \bullet \]
{\bf type D : } We have just one three-dimensional family of conifold singularities determined by the representation type $(1,(1);1,(1))$ so the fingerprint is $3_{con} \rTo \bullet$. As fingerprints are isomorphism invariants of the singularity, this finishes the proof.
\end{proof}

We claim that the minimal number of generators for these invariant rings is $7$. The structure of the invariant ring of three $2 \times 2$ matrices upto simultaneous conjugation was determined by Ed Formanek \cite{Formanek} who showed that it is generated by $10$ elements
\[
\{ tr(X_1),tr(X_2),tr(X_3),det(X_1),det(X_2),det(X_3),tr(X_1X_2),tr(X_1X_3),tr(X_2X_3),tr(X_1X_2X_3) \}
\]
and even gave the explicit quadratic polynomial satisfied by $tr(X_1X_2X_3)$ with coefficients in the remaining generators. The rings of invariants of the four cases of interest to us are quotients of this algebra by the ideal generated by three of its generators : for type $A$ it is $(det(X_1),det(X_2),det(X_3) )$, for type $B$ : $(det(X_1),tr(X_2),det(X_3) )$, for type $C$ :
$(det(X_1),tr(X_2),tr(X_3) )$ and for type $D$ : $( tr(X_1),tr(X_2),tr(X_3))$.

\section{Uniqueness of reduced setting} \label{unique}

In this section we will prove theorem~\ref{thmunique}. We will say that a vertex $v$ is {\em reducible} if one of the conditions $C^v_V$ (vertex removal), $C^v_l$ (loop removal in vertex dimension one) or $C^v_L$ (one (marked) loop removal) is satisfied. If we let the specific condition unspecified we will say that $v$ satisfies $C^v_X$ and denote $R^v_X$ for the corresponding marked quiver setting reduction. The resulting marked quiver setting will be denoted by
\[
R^v_X(Q^{\bullet},\alpha) \]
If $w \not= v$ is another vertex in $Q^{\bullet}$ we will denote the corresponding vertex in
$R^v_X(Q^{\bullet})$ also with $w$. The proof of the uniqueness result relies on three claims :
\begin{enumerate}
\item{If $w \not= v$ satisfies $R^w_Y$ in $(Q^{\bullet},\alpha)$, then $w$ virtually always satisfies
$R^w_Y$ in $R^v_X(Q^{\bullet},\alpha)$.}
\item{If $v$ satisfies $R^v_X$ and $w$ satisfies $R^w_Y$, then $R^v_X(R^w_Y(Q^{\bullet},\alpha)) = R^w_Y(R^v_X(Q^{\bullet},\alpha))$.}
\item{The previous two facts can be used to prove the result by induction on the minimal length of the reduction chain.}
\end{enumerate}
By the {\em neighborhood} of a vertex $v$ in $Q^{\bullet}$ we mean the (marked) subquiver on the
vertices connected to $v$. A neighborhood of a set of vertices is the union of the vertex-neighborhoods. {\em Incoming} resp. {\em outgoing} neighborhoods are defined in the natural manner.

\begin{lemma}\label{conditions}
Let $v \not= w$ be vertices in $(Q^{\bullet},\alpha)$. 
\begin{enumerate}
\item
If $v$ satisfies $C_V^v$ in $(Q^{\bullet},\alpha)$ and $w$ satisfies $C_{X}^w$,
then $v$ satisfies $C_{V}^w$ in $R_{X}^w(Q^{\bullet},\alpha)$ unless
the neighborhood of $\{ v,w \}$ looks like
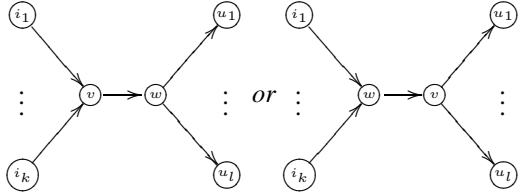
\[
\vcenter{\xymatrix@=.5cm{
\vtx{i_1}\ar[rd]&&&\vtx{u_1}\\
\vdots&\vtx{v}\ar[r]&\vtx{w}\ar[ru]\ar[rd]&\vdots\\
\vtx{i_k}\ar[ru]&&&\vtx{u_l}}
}\text{ or }
\vcenter{\xymatrix@=.5cm{
\vtx{i_1}\ar[rd]&&&\vtx{u_1}\\
\vdots&\vtx{w}\ar[r]&\vtx{v}\ar[ru]\ar[rd]&\vdots\\
\vtx{i_k}\ar[ru]&&&\vtx{u_l}}
}
\]
and $\alpha_v=\alpha_w$.
Observe that in this case $R_V^v(Q^{\bullet},\alpha) = R_V^w(Q^{\bullet},\alpha)$.
\item
If $v$ satisfies $C_{l}^v$ and $w$ satisfies $C_{X}^w$ then 
then $v$ satisfies $C_{l}^v$ in $R_{X}^w(Q^{\bullet},\alpha)$.
\item
If $v$ satisfies $C_{V}^v$ and $w$ satisfies $C_{X}^w$ then 
then $v$ satisfies $C_{V}^v$ in $R_{X}^w(Q^{\bullet},\alpha)$.
\end{enumerate} 
\end{lemma}

\begin{proof}
(1) : If $X=l$ then $R_{X}^w$ does not change the neighborhood of $v$ so $C_{V}^v$ holds in 
$R_{l}^w(Q^{\bullet},\alpha)$. If $X=L$ then $R_{X}^w$ does not change the neighborhood of $v$
unless $\alpha_v=1$ and $\chi_Q(\eps_w,\eps_v)=-1$ (resp. $\chi_Q(\eps_v,\eps_w)=-1$) depending on whether $w$ satisfies the in- or outgoing condition $C^w_L$. 
We only consider the first case, the latter is similar.
Then $v$ cannot satisfy the outgoing form of $C_V^v$ in $(Q^{\bullet},\alpha)$ so the incoming condition is satisfied. Because   
the $R_{L}^w$-move does not change the incoming neighborhood of $v$, $C_{V}^v$ still holds for $v$ in $R_{L}^w(Q^{\bullet},\alpha)$.

If $X=V$ and $v$ and $w$ have disjoint neighborhoods then $C_{V}^v$ trivially remains true in $R^w_V(Q^{\bullet},\alpha)$. Hence assume that there is at least one arrow from $v$ to $w$ (the case where there are only arrows from $w$ to $v$ is similar).
If $\alpha_v< \alpha_w$ then the incoming condition $C^v_V$ must hold (outgoing is impossible) and hence
$w$ does not appear in the incoming neighborhood of $v$. But then $R_V^w$ preserves the incoming neighborhood of $v$ and $C_V^v$ remains true in the reduction.    
If $\alpha_v> \alpha_w$ then the outgoing condition $C_V^w$ must hold and hence $w$ does not appear in the incoming neighborhood of $v$. So if the incoming condition $C_V^v$ holds  in $(Q^{\bullet},\alpha)$ it will still hold after the application of $R^w_{V}$.
If the outgoing condition $C^v_V$ holds, the neighborhoods of $v$ and $w$ in $(Q^{\bullet},\alpha)$ and
$v$ in $R_V^w(Q^{\bullet},\alpha)$ are depicted in figure~\ref{figure1001}
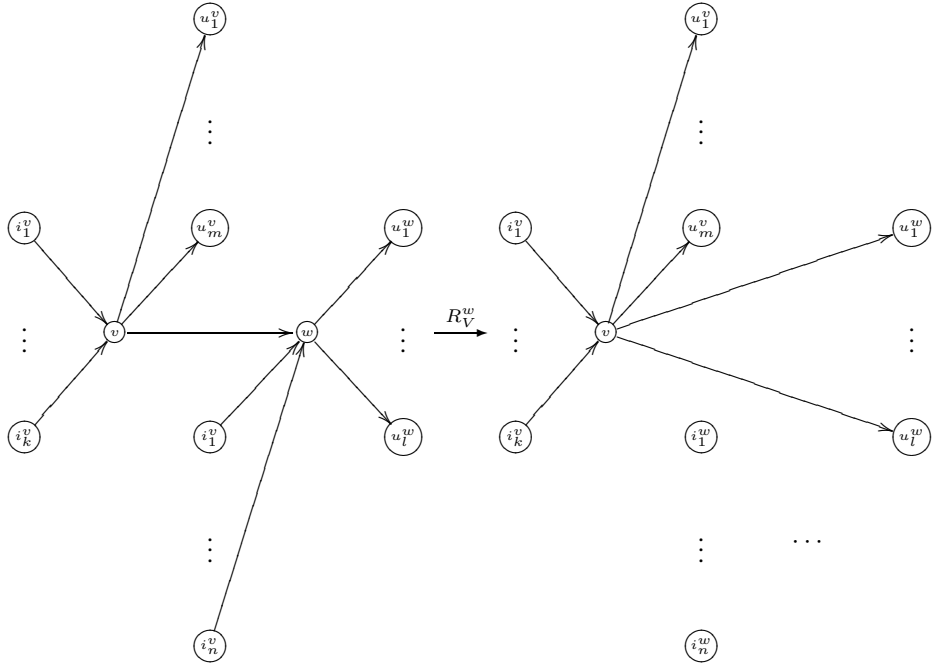
\begin{figure} 
\[
\vcenter{\xymatrix@=.8cm{
&&\vtx{u_1^v}&&\\
&&\vdots&&\\
\vtx{i_1^v}\ar[rd]&&\vtx{u_m^v}&&\vtx{u_1^w}\\
\vdots&\vtx{v}\ar[rr]\ar[ruuu]\ar[ru]&&\vtx{w}\ar[ru]\ar[rd]&\vdots\\
\vtx{i_k^v}\ar[ru]&&\vtx{i_1^v}\ar[ru]&&\vtx{u_l^w}\\
&&\vdots&&\\
&&\vtx{i_n^v}\ar[ruuu]&&
}} \rTo^{R^w_V} \vcenter{\xymatrix@=.8cm{
&&\vtx{u_1^v}&&\\
&&\vdots&&\\
\vtx{i_1^v}\ar[rd]&&\vtx{u_m^v}&&\vtx{u_1^w}\\
\vdots&\vtx{v}\ar[rrru]\ar[rrrd]\ar[ruuu]\ar[ru]&&&\vdots\\
\vtx{i_k^v}\ar[ru]&&\vtx{i_1^w}&&\vtx{u_l^w}\\
&&\vdots&\cdots&\\
&&\vtx{i_n^w}&&
}}
\]
\caption{Neighborhoods of $v$ and $w$}
\label{figure1001}
\end{figure}
Let $A$ be the set of arrows in $Q^{\bullet}$ and $A'$ the set of arrows in the reduction, then
because $\sum_{a \in A, s(a)=w}\alpha_{t(a)}\le \alpha_w$ (the incoming condition for $w$) we have
\begin{equation*}
\begin{split}
\sum_{a \in A', s(a)=v}\alpha'_{t(a)}&=
\sum_{a \in A,\atop s(a)=v, t(a)\ne w}\alpha_{t(a)}+\sum_{a \in A \atop t(a)=w,s(a)=v}\sum_{a \in A, s(a)=w}\alpha_{t(a)}\\
&\le 
\sum_{a \in A,\atop s(a)=v,t(a)\ne w}\alpha_{t(a)}+\sum_{a \in A \atop t(a)=w,s(a)=w}\alpha_w\\
&=\sum_{a \in A, s(a)=v}\alpha_{t(a)}\le \alpha_v
\end{split}
\end{equation*}
and therefore the outgoing condition $C^v_V$ also holds in $R_V^w(Q^{\bullet},\alpha)$.
Finally if $\alpha_v= \alpha_w$, it may be that 
$C_V^v$ does not hold  in $R_V^w(Q^{\bullet},\alpha)$. In this case $\chi(\epsilon_v,\alpha)<0$ and $\chi(\alpha,\epsilon_w)<0$ ($C_V^v$ is false in $R_V^w(Q^{\bullet},\alpha)$). 
Also $\chi(\alpha,\epsilon_v)\ge 0$ and $\chi(\epsilon_w,\alpha)\ge 0$ (otherwise $C_V$ does not hold for $v$ or $w$ in $(Q^{\bullet},\alpha)$). This implies that we are in the situation described in the lemma and the conclusion follows.

\noindent
(2) : None of the $R^w_X$-moves removes a loop in $v$ nor changes $\alpha_v = 1$.

\noindent
(3) : Assume that the incoming condition $C^v_L$ holds in $(Q^{\bullet},\alpha)$ but not in $R_X^w(Q^{\bullet},\alpha)$, then $w$ must be the unique vertex which has an arrow to $v$ and $X=V$. Because $\alpha_w=1<\alpha_v$,  the incoming condition $C^w_V$ holds. This means that there is also only one arrow arriving in $w$ and this arrow is coming from a vertex with dimension $1$. Therefore after applying  $R_V^w$, $v$ will still have only one incoming arrow starting in a vertex with dimension $1$.  A similar argument holds for the outgoing condition
$C^v_L$. 
\end{proof}

\begin{lemma}\label{moves}
Suppose that $v \not= w$ are vertices in $(Q^{\bullet},\alpha)$ and that $C_X^v$ and
$C_Y^w$ are satisfied.
If $C_X^v$ holds in $R^w_Y(Q^{\bullet},\alpha)$ and $C_Y^w$ holds in $R^v_X(Q^{\bullet},\alpha)$
then 
\[
R^v_X R^w_Y(Q^{\bullet},\alpha)= R^w_Y R^v_X(Q^{\bullet},\alpha)
\]
\end{lemma}

\begin{proof}
If $X,Y \in \{ l,L \}$ this is obvious, so let us assume that $X = V$. If $Y=V$ as well, we can calculate the Euler form $\chi_{R_V^w R_V^vQ}(\eps_{x},\eps_y)$.
Because
\[
\chi_{R_V^vQ}(\eps_{x},\eps_y)= \chi_{Q}(\eps_{x},\eps_y)- \chi_{Q}(\eps_{x},\eps_v)\chi_{Q}(\eps_v,\eps_y) 
\]
it follows that
\begin{equation*}
\begin{split}
\chi_{R_V^w R_V^vQ}(\eps_{x},\eps_y)&= \chi_{R_V^vQ}(\eps_{x},\eps_y)- \chi_{R_V^vQ}(\eps_{x},\eps_w)\chi_{R_V^vQ}(\eps_v,\eps_y)\\
&= \chi_{Q}(\eps_{x},\eps_y)- \chi_{Q}(\eps_{x},\eps_v)\chi_{Q}(\eps_v,\eps_y)\\& - \left(\chi_{Q}(\eps_{x},\eps_w)- \chi_{Q}(\eps_{x},\eps_v)\chi_{Q}(\eps_v,\eps_w)\right)
\left(\chi_{Q}(\eps_{w},\eps_y)- \chi_{Q}(\eps_{w},\eps_v)\chi_{Q}(\eps_v,\eps_y)\right)\\
&= \chi_{Q}(\eps_{x},\eps_y)- \chi_{Q}(\eps_{x},\eps_v)\chi_{Q}(\eps_v,\eps_y)- 
\chi_{Q}(\eps_{x},\eps_w)\chi_{Q}(\eps_{w},\eps_y) \\&- \chi_{Q}(\eps_{x},\eps_v)\chi_{Q}(\eps_v,\eps_w)\chi_{Q}(\eps_{w},\eps_v)\chi_{Q}(\eps_v,\eps_y)\\
& +\chi_{Q}(\eps_{x},\eps_w)\chi_{Q}(\eps_{w},\eps_v)\chi_{Q}(\eps_v,\eps_y)
+\chi_{Q}(\eps_{x},\eps_v)\chi_{Q}(\eps_v,\eps_w)\chi_{Q}(\eps_{w},\eps_y)
\end{split}
\end{equation*}
This is symmetric in $v$ and $w$ and therefore the ordering of $R_V^v$ and $R_V^w$ is irrelevant.

\noindent
If $Y=l$ we have the following equalities
\begin{equation*}
\begin{split}
\chi_{R_{l}^w R_V^vQ}(\eps_{x},\eps_y)&= \chi_{R_V^vQ}(\eps_{x},\eps_y)-\delta_{wx}\delta_{wy}\\
&=\chi_{Q}(\eps_{x},\eps_y)- \chi_{Q}(\eps_{x},\eps_v)\chi_{Q}(\eps_v,\eps_y)
-\delta_{wx}\delta_{wy}\\
&=\chi_{Q}(\eps_{x},\eps_y)-\delta_{wx}\delta_{wy}-(\chi_{Q}(\eps_{x},\eps_v)-\delta_{wx}\delta_{wv})(\chi_{Q}(\eps_v,\eps_y)-\delta_{wv}\delta_{wy})\\
&=\chi_{R_{l}^wQ}(\eps_{x},\eps_y) - \chi_{R_{l}^wQ}(\eps_{x},\eps_v)\chi_{R_{l}^wQ}(\eps_v,\eps_y)\\
&=\chi_{R_V^v R_{l}^wQ}.
\end{split}
\end{equation*}

\noindent
If $Y=L$, an $R_{L}^w$-move commutes with the $R_V^v$ move because it does not change the neighborhood of $v$ except when $v$ is the unique vertex of dimension $1$ connected to $w$. In this case
the neighborhood of $v$ looks like
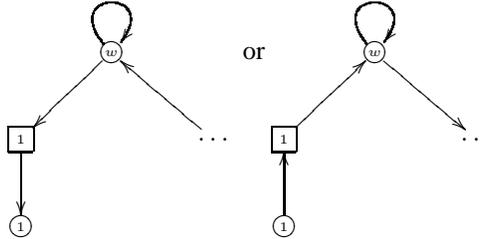
\vspace{.5cm}
\[
\xymatrix
{
&\vtx{w}\ar@(lu,ru)\ar[ld]&\\
\cvtx{1}\ar[d]&&\dots\ar[lu]\\
\vtx{1}&&}
\text{ or }
\xymatrix
{
&\vtx{w}\ar@(lu,ru)\ar[rd]&\\
\cvtx{1}\ar[ru]&&\dots\\
\vtx{1}\ar[u]&&}\]
In this case the reduction at $v$ is equivalent to a reduction at $v'$ (i.e. the lower vertex) which certainly commutes
with $R_{L}^w$. 
\end{proof}

\par \vskip 3mm
We are now in a position to prove theorem~\ref{thmunique}.

\begin{theorem}
If $(Q^{\bullet},\alpha)$ is a strongly connected marked quiver setting and 
$(Q_1^{\bullet},\alpha_1)$ and $(Q_2^{\bullet},\alpha_2)$ are two reduced marked quiver setting obtained
by applying reduction moves to $(Q^{\bullet},\alpha)$ then 
\[
(Q_1^{\bullet},\alpha_1) = (Q_2^{\bullet},\alpha_2)
\]
\end{theorem}

\begin{proof}
We do induction on the length $l_1$ of the reduction chain $R_1$ reducing $(Q^{\bullet},\alpha)$
to $(Q^{\bullet}_1,\alpha_1)$. If $l_1=0$, then $(Q^{\bullet},\alpha)$ has no reducible vertices so the result holds trivially. Assume the result holds for all lengths $< l_1$. There are two cases to consider.

There exists a vertex $v$ satisfying a loop removal condition $C_{X}^v, X=l$ or $L$. Then,
there is a $R_X^v$-move in both reduction chains $R_1$ and $R_2$. This follows from lemma \ref{conditions}
and the fact that none of the vertices in $(Q_1^{\bullet},\alpha_1)$ and $(Q_2^{\bullet},\alpha_2)$ are reducible. By the commutation relations from lemma \ref{moves}, we can bring this reduction to the first position in both chains and use induction.

If there is a vertex $v$ satisfying condition $C_V^v$, either both chains will contain 
an $R_V^v$-move or the neighborhood of $v$ looks like the figure in lemma \ref{conditions} (1).
Then, $R_1$ can contain an $R_V^v$-move and $R_2$ an $R_V^w$-move. But then we change the $R_V^w$ move into a $R_V^v$ move, because they have the same effect. The
concluding argument is similar to that above.
\end{proof}

\section{The main result}

\begin{theorem} There are no reduced quiver singularities for $d \leq 2$. For $d=3$ the conifold
is the only reduced quiver singularity. For $d=4$ (resp. $d=5$ and $d=6$) there are precisely
three (resp. ten and fifty-three) reduced quiver singularities.
\end{theorem}

\begin{proof}
The details for $d \leq 4$ were given above. The classification for $d=5$ and $d=6$ is given in the full version of this paper which is available at 

{\tt ftp://wins.uia.ac.be/pub/preprints/02/SOSfull.pdf}
\end{proof}

\begingroup\raggedright\endgroup


\begin{thebibliography}{10}

\bibitem{Bocklandt}
Raf Bocklandt, {\it Smooth quiver quotient varieties}, {\em Linear Alg. Appl.}
(2002) (to appear)  {\tt
  arXiv:math.RT/0204355}.

\bibitem{BockSym}
Raf Bocklandt and Stijn Symens, {\it Isolated quiver singularities}, (in preparation)

\bibitem{Formanek}
Edward Formanek, {\it Invariants and the ring of generic matrices}, {\em J. Algebra}
{\bf 89} (1984) 178-223


\bibitem{LBsmooth}
Lieven Le Bruyn, {\it Local structure of Schelter-Procesi smooth orders}, {\em Trans. AMS}
{\bf 352} (2000) 4815-4841

\bibitem{LBProcesi}
Lieven Le Bruyn and Claudio Procesi, {\it Semi-simple representations of quivers}, {\em Trans. AMS} {\bf 317} (1990) 585-598

\bibitem{PORTA} PORTA - A Polyhedron Representation and Transformation 
               Algorithm {\tt
http://elib.zib.de/pub/Packages/mathprog/polyth/porta/}

\bibitem{ProcesiCH}
Claudio Procesi, {\it A formal inverse to the Cayley-Hamilton theorem}, {\em J. Algebra}
{\bf 107} (1987) 143-176

\bibitem{Schelter}
William Schelter, {\it Smooth algebras}, {\em J. Algebra} {\bf 103} (1986) 677-685





\end{thebibliography}
\end{document}